\input amstex.tex
\documentstyle{amsppt}
\magnification=\magstep 1
\TagsAsMath

\define\flow{\left(\bold{M},\{S^t\}_{t\in\Bbb R},\mu\right)}
\define\symb{\Sigma=\left(\sigma(1),\dots,\sigma(m)\right)}
\define\jflow{\left(\bold{M}_j,\{S_j^t\},\mu_j\right)}

\document

{\catcode`\@=11\gdef\logo@{}}

\noindent
June 21, 1999

\bigskip \bigskip

\centerline{\bf The Complete Hyperbolicity of Cylindric Billiards}

\bigskip \bigskip \bigskip

\centerline{{\bf N\'andor Sim\'anyi}
\footnote{Research supported by the Hungarian National Foundation for
Scientific Research, grants OTKA-26176 and OTKA-29849.}}
\centerline{Bolyai Institute, University of Szeged,}
\centerline{6720 Szeged, Aradi V\'ertanuk tere 1, Hungary.}
\centerline{E-mail: simanyi\@math.u-szeged.hu}

\bigskip \bigskip

\hbox{\centerline{\vbox{\hsize 8cm {\bf Abstract.}
The connected configuration space of a so called cylindric
billiard system is a flat torus minus finitely many spherical
cylinders. The dynamical system describes the uniform motion of a point
particle in this configuration space with specular reflections at the
boundaries of the removed cylinders. It is proven here that under a
certain geometric condition --- slightly stronger than the necessary
condition presented in [S-Sz(1998)] --- a cylindric billiard flow is
completely hyperbolic. As a consequence, every hard ball system is completely
hyperbolic --- a result strengthening the theorem of [S-Sz(1999)].}}}

\bigskip \bigskip

\heading
1. Introduction
\endheading

\bigskip \bigskip

Non-uniformly hyperbolic systems (possibly, with singularities) play a pivotal
role in the ergodic theory of dynamical systems. Their systematic study
started several decades ago, and it is not our goal here to provide the reader
with a comprehensive review of the history of these investigations but,
instead, we opt for presenting in nutshell a cross section of a few selected
results.

In 1939 G. A. Hedlund and E. Hopf [He(1939)], [Ho(1939)], proved the
hyperbolic ergodicity of geodesic flows on closed, compact surfaces with
constant negative curvature by inventing the famous method of "Hopf chains"
constituted by local stable and unstable invariant manifolds.

In 1963 Ya. G. Sinai [Sin(1963)] formulated a modern version of Boltzmann's
ergodic hypothesis, what we call now the "Boltzmann-Sinai ergodic hypothesis":
the billiard system of $N$ ($\ge2$) hard balls of unit mass moving in the
flat torus $\Bbb T^\nu=\Bbb R^\nu/\Bbb Z^\nu$ ($\nu\ge2$) is ergodic after
we make the standard reductions by fixing the values of the trivial invariant
quantities. It took seven years until he proved this conjecture for the
case $N=2$, $\nu=2$ in [Sin(1970)]. Another 17 years later N. I. Chernov
and Ya. G. Sinai [S-Ch(1987)] proved the hypothesis for the case $N=2$,
$\nu\ge2$ by also proving a powerful and very useful theorem on local
ergodicity.

In the meantime, in 1977, Ya. Pesin [P(1977)] laid down the foundations
of his theory on the ergodic properties of smooth, hyperbolic dynamical
systems. Later on this theory (nowadays called Pesin theory) was
significantly extended by A. Katok and J-M. Strelcyn [K-S(1986)]
to hyperbolic systems with singularities. That theory is already
applicable for billiard systems, too.

Until the end of the seventies the phenomenon of hyperbolicity (exponential
unstability of the trajectories) was almost exclusively attributed to some
direct geometric scattering effect, like negative curvature of space, or
strict convexity of the scatterers. This explains the profound shock that
was caused by the discovery of L. A. Bunimovich [B(1979)]: certain focusing
billiard tables (like the celebrated stadium) can also produce complete
hyperbolicity and, in that way, ergodicity. It was partly this result that
led to Wojtkowski's theory of invariant cone fields, [W(1985)], 
[W(1986)].

The big difference between the system of two balls in $\Bbb T^\nu$
($\nu\ge2$, [S-Ch(1987)]) and the system of $N$ ($\ge3$) balls in
$\Bbb T^\nu$ is that the latter one is merely a so called semi-dispersive
billiard system (the scatterers are convex but not strictly convex
sets, namely cylinders), while the former one is strictly dispersive
(the scatterers are strictly convex sets). This fact makes the proof
of ergodicity (mixing properties) much more complicated. In our series
of papers jointly written with A. Kr\'amli and D. Sz\'asz [K-S-Sz(1990)],
[K-S-Sz(1991)], and [K-S-Sz(1992)] we managed to prove the (hyperbolic)
ergodicity of three and four billiard balls in the toroidal container
$\Bbb T^\nu$. By inventing new topological methods and the Connecting Path
Formula (CPF), in my two-part paper [Sim(1992)] I proved the (hyperbolic)
ergodicity of $N$ hard balls in $\Bbb T^\nu$, provided that $N\le\nu$.

The common feature of hard ball systems is --- as D. Sz\'asz pointed this
out first in [Sz(1993)] and [Sz(1994)] --- that all of theom belong to the
family of so called cylindric billiards, the definition of which can be
found later in this paragraph. However, the first appearance of a special,
3-D cylindric billiard system took place in [K-S-Sz(1989)], where we
proved the ergodicity of a 3-D billiard flow with two orthogonal
cylindric scatterers. Later D. Sz\'asz [Sz(1994)] presented a complete
picture (as far as ergodicity is concerned) of cylindric billiards with
cylinders whose generator subspaces are spanned by mutually orthogonal
coordinate axes. The task of proving ergodicity for the first non-trivial,
non-orthogonal cylindric billiard system was taken up in [S-Sz(1994)].

Finally, in our joint venture with D. Sz\'asz [S-Sz(1999)] we managed to
prove the complete hyperbolicity of {\it typical} hard ball systems.

\subheading{\bf 1.1. Cylindric billiards} Consider the $d$-dimensional
($d\ge2$) flat torus $\Bbb T^d=\Bbb R^d/\Cal L$ supplied with the
usual Riemannian inner product $\langle\, .\, ,\, .\, \rangle$ inherited
from the standard inner product of the universal covering space $\Bbb R^d$.
Here $\Cal L\subset\Bbb R^d$ is supposed to be a lattice, i. e. a discrete
subgroup of the additive group $\Bbb R^d$ with $\text{rank}(\Cal L)=d$.
The reason why we want to allow general lattices other than just the
integer lattice $\Bbb Z^d$ is that otherwise the hard ball systems would
not be covered! The geometry of the structure lattice $\Cal L$ in the
case of a hard ball system is significantly different from the geometry
of the standard lattice $\Bbb Z^d$ in the standard Euclidean space
$\Bbb R^d$, see subsection 2.4, especially (2.4.2) and (2.4.5).

The configuration space of a cylindric billiard is
$\bold Q=\Bbb T^d\setminus\left(C_1\cup\dots\cup C_k\right)$, where the
cylindric scatterers $C_i$ ($i=1,\dots,k$) are defined as follows:

Let $A_i\subset\Bbb R^d$ be a so called lattice subspace of $\Bbb R^d$,
which means that $\text{rank}(A_i\cap\Cal L)=\text{dim}A_i$. In this case
the factor $A_i/(A_i\cap\Cal L)$ is a subtorus in $\Bbb T^d=\Bbb R^d/\Cal L$
which will be taken as the generator of the cylinder 
$C_i\subset\Bbb T^d$, $i=1,\dots,k$. Denote by $L_i=A_i^\perp$ the
orthocomplement of $A_i$ in $\Bbb R^d$. Throughout this article we will
always assume that $\text{dim}L_i\ge2$. Let, furthermore, the numbers
$r_i>0$ (the radii of the spherical cylinders $C_i$) and some translation
vectors $t_i\in\Bbb T^d=\Bbb R^d/\Cal L$ be given. The translation
vectors $t_i$ play a crucial role in positioning the cylinders $C_i$
in the ambient torus $\Bbb T^d$. Set

$$
C_i=\left\{x\in\Bbb T^d:\; \text{dist}\left(x-t_i,A_i/(A_i\cap\Cal L)
\right)<r_i \right\}.
$$
In order to avoid further unnecessary complications, we always assume that
the interior of the configuration space 
$\bold Q=\Bbb T^d\setminus\left(C_1\cup\dots\cup C_k\right)$ is connected.
The phase space $\bold M$ of our cylindric billiard flow will be the
unit tangent bundle of $\bold Q$ (modulo some natural glueings at its
boundary), i. e. $\bold M=\bold Q\times\Bbb S^{d-1}$. (Here $\Bbb S^{d-1}$
denotes the unit sphere of $\Bbb R^d$.)

The dynamical system $\flow$, where $S^t$ ($t\in\Bbb R$) is the dynamics 
defined by uniform motion inside the domain $\bold Q$ and specular
reflections at its boundary (at the scatterers), and $\mu$ is the
Liouville measure, is called a cylindric billiard flow we want to
investigate. (As to notions and notations in connection with 
semi-dispersive billiards, the reader is kindly recommended to consult
the work [K-S-Sz(1990)].)

\subheading{\bf Transitive cylindric billiards}

The main conjecture concerning the (hyperbolic) ergodicity of cylindric
billiards is the "Erd\H otarcsa conjecture" (named after the picturesque
village in rural Hungary where it was initially formulated) that appeared
as Conjecture 1 in Section 3 of [S-Sz(1998)]:

\subheading{\bf The Erd\H otarcsa conjecture} A cylindric billiard flow is
ergodic if and only if it is transitive. (As for the definition and basic
features of transitivity, see Section 3 (especially between 3.1 and 3.6)
of [S-Sz(1998)] or subsection 2.2 below.) In that case the cylindric billiard
system is actually a completely hyperbolic Bernoulli flow, see
[C-H(1996)] and [O-W(1998)].

The theorem of this paper proves a slightly relaxed version of this
conjecture (only full hyperbolicity without ergodicity) for a wide class
of cylindric billiard systems, namely the so called "transverse systems"
(see subsection 2.3 below) which include every hard ball system:

\subheading{\bf Theorem} Assume that the cylindric billiard system is
transverse, see subsection 2.3. Then this billiard flow is completely
hyperbolic, i. e. all relevant Lyapunov exponents are nonzero almost
everywhere. Consequently, such dynamical systems have (at most countably
many) ergodic components of positive measure, and the restriction of
the flow to the ergodic components has the Bernoulli property,
see [C-H(1996)] and [O-W(1998)].

\subheading{\bf Corollary of the theorem} Every hard ball system ---
necessarily being a transverse cylindric billiard system, see 
subsection 2.4 --- is completely hyperbolic.

Thus, the theorem of this paper generalizes the main result of 
[S-Sz(1999)], where the complete hyperbolicity of {\it almost every}
hard ball system was proven.

\subheading{\bf Organizing of the paper} After the technical preparation in
Section 2, the theorem will be proven in the two subsequent sections.
According to the usually accepted strategy developed in the series of papers
[K-S-Sz(1989, 1991, 1992)], [Sim(1992)], the proof of full hyperbolicity
should consist of two major steps:

\subheading{\bf Step 1. (Geometric-algebraic considerations)} To prove
that the existence of a combinatorially rich (appropriately defined!)
trajectory segment $S^{[a,b]}x_0$ for a smooth phase point $x_0\in\bold M$
implies (modulo some smooth, proper submanifolds) that the phase point
$x_0$ is hyperbolic (or, using the older language, sufficient). This will
be carried out in Section 3.

\subheading{\bf Step 2. (Dynamical-topological part)} To show that for
$\mu$-almost every phase point $x_0$ the symbolic collision sequence of
the entire trajectory of $x_0$ is combinatorially rich. This will be
accomplished in Section 4.

\bigskip \bigskip

\heading
2. Prerequisites
\endheading

\bigskip \bigskip

\subheading{\bf 2.1. Sub-billiards} Assume that a subset 
$I\subset\{1,\dots,k\}$ is given, and we consider the cylindric billiard flow
in the torus $\Bbb T^d=\Bbb R^d/\Cal L$ so that only the cylinders
$\{C_i:\, i\in I\}$ are retained as scatterers; the other ones are no longer
removed from the configuration space $\bold Q$ and the uniformly moving point
$q=q(x)$ ($x=(q,v)\in\bold M$) can freely pass through them. We call the
arising billiard flow a {\it sub-billiard}.

It turns out pretty soon that the name "factor billiard" would have been much
better. Namely, let us consider the linear subspaces 
$E^+=\text{span}\{L_i:\, i\in I\}$ and 
$E_0=\left(E^+\right)^\perp=\bigcap_{i\in I}A_i$. It is an elementary exercise
to show that the intersection $E_0$ of the lattice subspaces $A_i$ is also a
lattice subspace, i. e. $\text{rank}(E_0\cap\Cal L)=\text{dim}E_0$.

It is easy to see that the sub-billiard flow $\left\{S_I^t\right\}$
($t\in\Bbb R$) defined by the scatterers $\{C_i:\, i\in I\}$ has the
following peculiarity: the velocity component $P_{E_0}(v_t)$ of the moving
phase point $x_t=(q_t,v_t)$ does not change, and in the direction of the
subspace $E_0$ (or, equivalently, in the direction of the subtorus 
$E_0/(E_0\cap\Cal L)\subset\Bbb R^d/\Cal L$) the motion of $q_t$ is
conditionally periodic. (Here, as always, $P_{E_0}(\, .\, )$ denotes
the orthogonal projection of $\Bbb R^d$ onto the subspace $E_0$.)
According to the invariance of the quantity $P_{E_0}(v_t)$, we fix its
value by introducing the reduction $P_{E_0}(v_t)=0$. After this reduction
the sub-billiard flow $\left\{S_I^t\right\}$ will have a translation
invariance in the direction of the subtorus $E_0/(E_0\cap\Cal L)$, thus we
factorize out the configuration space with respect to spatial translations
by elements $\tau\in E_0/(E_0\cap\Cal L)$ as follows:
$q\sim q'\Longleftrightarrow q-q'\in E_0/(E_0\cap\Cal L)$. The flow arising
after the reduction $P_{E_0}(v)=0$ and the above factorization is denoted by
$\left\{S_I^t\right\}$. Let us describe now the natural configuration and
velocity spaces of the flow $\left\{S_I^t\right\}$. The velocity space is
obviously the orthocomplement $E^+$ of the lattice subspace 
$E_0\subset\Bbb R^d$. (We note that the space $E^+$ does not have to be a
lattice subspace!) After specifying the kinetic energy
$\varepsilon=\frac{1}{2}||v||^2$ of the subsystem, we get the sphere of
radius $\sqrt{2\varepsilon}$ in the Euclidean space $E^+$ as the velocity
space for the sub-billiard flow $\left\{S_I^t\right\}$. As far as the
configuration space $\bold Q=\bold Q_I$ is concerned, it is naturally the
factor torus

$$
\Bbb T^d/\left(E_0/(E_0\cap\Cal L)\right)=\Bbb R^d/(\Cal L+E_0)
$$
(minus the intersections of the cylinders $\{C_i:\, i\in I\}$ with that
factor torus) supplied with the Euclidean metric of the space $E^+$
as the Riemannian metric on $\Bbb T^d/\left(E_0/(E_0\cap\Cal L)\right)$.
Note that the subspace $E^+$ can be naturally identified with the tangent
spaces of the factor torus $\Bbb T^d/\left(E_0/(E_0\cap\Cal L)\right)$
at different points. 

By projecting the whole space $\Bbb R^d$ onto $E^+$
we see that the factor torus 

$$
\Bbb T^d/\left(E_0/(E_0\cap\Cal L)\right)=\Bbb R^d/(\Cal L+E_0)
$$
can be naturally identified with the factor $E^+/P_{E^+}(\Cal L)$. We note
that --- as it follows easily from the fact that $E_0=(E^+)^\perp$ is a
lattice subspace --- the projection $P_{E^+}(\Cal L)$ of the lattice
$\Cal L$ onto $E^+$ is a lattice in the subspace $E^+$.

\subheading{\bf 2.2. Transitivity} Let $L_1,\dots,L_k\subset\Bbb R^d$
be subspaces, $\text{dim}L_i\ge2$, $A_i=L_i^\perp$, $i=1,\dots,k$. Set

$$
\Cal G_i=\left\{U\in\text{SO}(d):\, U|A_i=\text{Id}_{A_i}\right\},
$$
and let 
$\Cal G=\left\langle\Cal G_1,\dots,\Cal G_k\right\rangle\subset\text{SO}(d)$
be the algebraic generate of the compact, connected Lie subgroups
$\Cal G_i$ in $\text{SO}(d)$. The following notions appeared in Section 3 of
[S-Sz(1998)].

\subheading{\bf Definition 2.2.1} We say that the system of base spaces
$\{L_1,\dots,L_k\}$ (or, equivalently, the cylindric billiard system defined
by them) is {\it transitive} if and only if the group $\Cal G$ acts 
transitively on the unit sphere $\Bbb S^{d-1}$ of $\Bbb R^d$.

\subheading{Definition 2.2.2} We say that the system of subspaces
$\{L_1,\dots,L_k\}$ has the Orthogonal Non-splitting Property (ONSP) if there
is no non-trivial orthogonal splitting $\Bbb R^d=B_1\oplus B_2$ of
$\Bbb R^d$ with the property that for every index $i$ ($1\le i\le k$)
$L_i\subset B_1$ or $L_i\subset B_2$.

The next result can be found in Section 3 of [S-Sz(1998)] (see 3.1--3.6
thereof):

\subheading{\bf Proposition 2.2.3} For the system of subspaces
$\{L_1,\dots,L_k\}$ the following three properties are equivalent:

(1) $\{L_1,\dots,L_k\}$ is transitive;

(2) $\{L_1,\dots,L_k\}$ has the ONSP;

(3) the natural representation of $\Cal G$ in $\Bbb R^d$ is irreducible.

\subheading{\bf 2.3. Transverseness}

\subheading{\bf Definition 2.3.1} We say that the system of subspaces
$\{L_1,\dots,L_k\}$ of $\Bbb R^d$ is {\it transverse} if the following
property holds: For every {\it non-transitive} subsystem
$\{L_i:\, i\in I\}$ ($I\subset\{1,\dots,k\}$) there exists an index 
$j_0\in\{1,\dots,k\}$ such that $P_{E^+}(A_{j_0})=E^+$, where
$A_{j_0}=L^\perp_{j_0}$, and $E^+=\text{span}\{L_i:\, i\in I\}$. We note
that in this case, necessarily, $j_0\not\in I$, otherwise $P_{E^+}(A_{j_0})$
would be orthogonal to the subspace $L_{j_0}\subset E^+$. Therefore, every
transverse system is automatically transitive.

\subheading{\bf 2.4. A major family of examples}

\subheading{\bf 2.4.1. Hard ball systems} Hard ball systems in the standard
unit torus $\Bbb T^\nu=\Bbb R^\nu/\Bbb Z^\nu$ ($\nu\ge2$) with positive masses
$m_1,\dots,m_N$ are described (for example) in Section 1 of [S-Sz(1999)].
These are the dynamical systems describing the motion of $N$ ($\ge2$) hard
balls with radius $r>0$ and positive masses $m_1,\dots,m_N$ in the 
standard unit torus $\Bbb T^\nu=\Bbb R^\nu/\Bbb Z^\nu$. The center of the
$i$-th ball is denoted by $q_i$ ($\in\Bbb T^\nu$), its time derivative is
$v_i=\dot q_i$, $i=1,\dots,N$. One uses the standard reduction of kinetic
energy $\varepsilon=\frac{1}{2}\sum_{i=1}^N m_i||v_i||^2=\frac{1}{2}$.
The arising configuration space (still without the removal of the scattering
cylinders $C_{i,j}$) is the torus

$$
\Bbb T^{\nu N}=\left(\Bbb T^{\nu}\right)^N=\left\{(q_1,\dots,q_N):\;
q_i\in\Bbb T^\nu,\; i=1,\dots,N\right\}
$$
supplied with the Riemannian inner product 

$$
\langle v,v'\rangle=\sum_{i=1}^N m_i\langle v_i,v'_i \rangle
\tag 2.4.2
$$
in its common tangent space $\Bbb R^{\nu N}=\left(\Bbb R^{\nu}\right)^N$.
Now the Euclidean space $\Bbb R^{\nu N}$ with the inner product (2.4.2)
plays the role of $\Bbb R^d$ in the original definition of cylindric
billiards, see Section 1 above.

The generator subspace $A_{i,j}\subset \Bbb R^{\nu N}$ ($1\le i<j\le N$)
of the cylinder $C_{i,j}$ (describing the collisions between the $i$-th and
$j$-th balls) is given by the equation

$$
A_{i,j}=\left\{(q_1,\dots,q_N)\in\left(\Bbb R^\nu\right)^N:\; 
q_i=q_j \right\},
\tag 2.4.3
$$
see (4.3) in [S-Sz(1998)].
Its orthocomplement $L_{i,j}\subset\Bbb R^{\nu N}$ is then defined by the
equation

$$
L_{i,j}=\left\{(q_1,\dots,q_N)\in\left(\Bbb R^\nu\right)^N:\; 
q_k=0 \text{ for } k=i,j, \text{ and } m_iq_i+m_jq_j=0 \right\},
\tag 2.4.4
$$
see (4.4) in [S-Sz(1998)].
Easy calculation shows that the cylinder $C_{i,j}$ is indeed spherical and
the radius of its base sphere is equal to
$r_{i,j}=2r\sqrt{\frac{m_im_j}{m_i+m_j}}$, see Section 4, especially formula
(4.6) in [S-Sz(1998)].

The structure lattice $\Cal L\subset\Bbb R^{\nu N}$ is clearly the integer
lattice $\Cal L=\Bbb Z^{\nu N}$. 

Due to the presence of an extra invariant quantity
$I=\sum_{i=1}^N m_iv_i$, one usually makes the reduction
$\sum_{i=1}^N m_iv_i=0$ and, correspondingly, factorizes the configuration
space with respect to uniform spatial translations:

$$
(q_1,\dots,q_N)\sim(q_1+a,\dots,q_N+a), \quad a\in\Bbb T^\nu,
$$
see also subsection 2.1 above. The natural, common tangent space of this
reduced configuration space is then

$$
\Cal Z=\left\{(v_1,\dots,v_N)\in\left(\Bbb R^\nu\right)^N:\;
\sum_{i=1}^N m_iv_i=0\right\}=\left(\bigcap_{i<j}A_{i,j}
\right)^\perp=\left(\Cal A\right)^\perp
\tag 2.4.5
$$
supplied again with the inner product (2.4.2), see also (4.1) and
(4.2) in [S-Sz(1998)]. The base spaces $L_{i,j}$
of (2.4.4) are obviously subspaces of $\Cal Z$, and we take
$\tilde A_{i,j}=A_{i,j}\cap\Cal Z=P_{\Cal Z}(A_{i,j})$ as the orthocomplement
of $L_{i,j}$ in $\Cal Z$.

Note that the configuration space of the reduced system (with 
$\sum_{i=1}^N m_iv_i=0$) is naturally the torus
$\Bbb R^{\nu N}/(\Cal A+\Bbb Z^{\nu N})=\Cal Z/P_{\Cal Z}(\Bbb Z^{\nu N})$,
see also subsection 2.1.

\subheading{Proposition 2.4.6} For every hard ball system with parameters
$N,\nu,r,m_1,\dots\allowmathbreak,m_N$
($N,\nu\ge2$, $m_i,r>0$) the collection of base
spaces $\{L_{i,j}:\, 1\le i<j\le N\}$ has the property of transverseness
in the tangent space $\Cal Z$.

\subheading{\bf Proof} Assume that 
$I\subset\left\{(i,j):\, 1\le i<j\le N\right\}$ is the index set of a
non\-tran\-si\-tive family $\{L_{i,j}:\, (i,j)\in I\}$ of base spaces in 
$\Cal Z$. The set $I$ can be considered as the set of edges of a
non-oriented collision graph $\Cal G$ with vertex set $\{1,\dots,N\}$.
It is shown in Remark 4.12 of [S-Sz(1998)] that the non-transitivity
of $\{L_{i,j}:\, (i,j)\in I\}$ means that the graph $\Cal G$ is not
connected on the full vertex set $\{1,\dots,N\}$. Choose a pair $(i_0,j_0)$
($1\le i_0<j_0\le N$) so that these indices belong to different connected
components of $\Cal G$. Then elementary consideration shows that
$P_{E^+}(\tilde A_{i_0,j_0})=E^+$, where
$E^+=\text{span}\{L_{i,j}:\, (i,j)\in I\}$.
(As a matter of fact, specifying an element $q\in E^+$ means specifying
the relative positions of the balls in each connected component of the
graph $\Cal G$. Finding a suitable element
$\tilde q\in\tilde A_{i_0,j_0}$ with $P_{E^+}(\tilde q)=q$ precisely means
that we ought to move each of the connected components of $\Cal G$
uniformly in the ambient torus so that the centers of the $i_0$-th and 
$j_0$-th balls just coincide. However, this can obviously be accomplished.)

This finishes the proof of the proposition. \qed

\subheading
{\bf 2.5. Another family of examples: Connected "direct sum systems"} 

Consider now such cylindric billiard systems in which the space
$\Bbb R^d$ decomposes into a linear direct sum

$$
\Bbb R^d=L_1+L_2+\dots+L_k
\tag 2.5.1
$$
of the base spaces $L_i$. With the decomposition (2.5.1) we associate a
non-oriented graph $\Cal G$ with the vertex set
$\Cal V(\Cal G)=\{1,\dots,k\}$ and edge set

$$
\Cal E(\Cal G)=\left\{\{i,j\}:\, i\not=j\text{ and }L_i
\not\perp L_j\right\}.
$$
It is then obvious that the transitivity of such a cylindric billiard system
is equivalent to the connectedness of the graph of non-orthogonality
$\Cal G$ (on the full vertex set $\{1,\dots,k\}$), which we assume now.

\subheading{\bf Proposition 2.5.2} A "direct sum system" (described above)
with a connected graph of non-orthogonality $\Cal G$ enjoys the property
of transverseness.

\subheading{Proof} Assume that $I\subset\{1,\dots,k\}$, and the system of
subspaces $\{L_i:\, i\in I\}$ is not transitive, i. e. $|I|<k$. Now,
for any index $j_0\in\{1,\dots,k\}\setminus I$ one has

$$
L_{j_0}\cap\text{span}\left\{L_i:\; i\in I\right\}=L_{j_0}\cap E^+=
\{0\},
$$
i. e. $\text{span}\left\{A_{j_0},(E^+)^\perp\right\}=\Bbb R^d$ which, in turn,
means that $P_{E^+}(A_{j_0})=E^+$. \qed

\subheading{\bf Remark} Consider a hard ball system with the graph of allowed
collisions $\Cal G$, see Remark 4.12 in [S-Sz(1998)]. Assume that the graph
$\Cal G$ is a tree, i. e. a connected graph without loop. It is an easy
exercise to see that such a hard ball system belongs to the family of
connected direct sum systems described above.

\subheading{\bf Hyperbolic (sufficient) trajectories} Their definition and
fundamental properties can be found --- for example --- in Definition 2.12
and Lemma 2.13 of [K-S-Sz(1990)].

\subheading{\bf 2.7. The subsets $\bold M^0$ and $\bold M^\#$} Denote by
$\bold M^\#$ the set of all phase points $x\in\bold M$ for which the
trajectory of $x$ encounters infintely many non-tangential collisions
in both time directions. The trajectories of the points 
$x\in\bold M\setminus\bold M^\#$ are lines: the motion is linear and uniform,
see the appendix of [Sz(1994)]. It is proven in lemmas A.2.1 and A.2.2
of [Sz(1994)] that the closed set $\bold M\setminus\bold M^\#$ is a finite
union of hyperplanes. Thus, in our study of complete hyperbolicity,
we can discard the set $\bold M\setminus\bold M^\#$ and focus
on the open set $\bold M^\#$ with full measure.

Denote by $\bold M^0$ the set of all non-singular phase points $x\in\bold M$,
i. e. all phase points $x$ whose entire trajectory is smooth. Since the
complement $\bold M\setminus\bold M^0$ of this set is a countable union of
smooth, proper submanifolds of $\bold M$, we can again discard the zero set
$\bold M\setminus\bold M^0$ and only consider phase points
$x\in\bold M^0\cap\bold M^\#$.

\subheading{\bf Finitely many collisions in finite time}
By the results of Vaserstein [V(1979)], Galperin [G(1981)] and
Burago-Ferleger-Kononenko [B-F-K(1998)], in a 
semi-dis\-per\-sive billiard flow there can only be finitely many 
collisions in finite time intervals, see Theorem 1.1 in [B-F-K(1998)]. 
Thus, the dynamics is well defined as long as the trajectory does not hit 
more than one boundary components at the same time.

\bigskip \bigskip

\heading
3. Geometric-Algebraic Considerations.
\endheading

\bigskip \bigskip

We begin this section with some new notions. Consider the linear subspaces
$A_i$ and $L_i=A^\perp_i$ ($i=1,\dots,k$, $\text{dim}L_i\ge 2$) in $\Bbb R^d$
and the positive numbers (radii) $r_i$ associated with them. Furthermore,
consider and fix a finite sequence $\Sigma=(\sigma(1),\dots,\sigma(m))$
of labels $\sigma(j)\in\{1,2,\dots,k\}$, a so called symbolic collision
sequence.

\subheading{\bf Definition 3.1} We say that $\gamma$ is a {\it Euclidean path}
with the collision sequence $\Sigma$ if the following properties hold:

(1) $\gamma:\, [0,\infty)\to \Bbb R^d$ is a piecewise linear, continuous curve
in $\Bbb R^d$ with $\gamma(0)=0$;

(2) $\gamma$ has an arc length parametrization by $t$, i. e. 
$||\dot\gamma(t)||=1$ for $t\ge 0$;

(3) the velocity $\dot\gamma(t)$ has finitely many discontinuities and these
discontinuities are jump discontinuities taking place at time moments
($0<$)$t_1<t_2<\dots<t_m<\infty$;

(4) the vectors of abrupt velocity change 
$\dot\gamma(t_j+0)-\dot\gamma(t_j-0)\ne 0$ belong to the base subspace
$L_{\sigma(j)}$, $j=1,\dots,m$.

We can think of the curve $\gamma$ as the Euclidean lifting of a finite
trajectory segment (extended to $t\to\infty$ with constant velocity, just
for technical reasons) of the genuine cylindric billiard flow in
$\Bbb T^d=\Bbb R^d/\Cal L$. The $j$-th collision takes place at time moment
$t_j$ at the boundary of the translated cylinder 
$a_j+C_{\sigma(j)}=a_j(\gamma)+C_{\sigma(j)}$ ($j=1,\dots,m$), where 

$$
\aligned
C_i=\left\{x\in\Bbb R^d:\, \text{dist}(x,A_i)<r_i\right\} \; \;
(i=1,\dots,k), \\
a_j=a_j(\gamma)=P_{\sigma(j)}\left(\gamma(t_j)\right)-r_{\sigma(j)}\cdot
\frac{\dot\gamma(t_j+0)-\dot\gamma(t_j-0)}
{\left\Vert\dot\gamma(t_j+0)-\dot\gamma(t_j-0)\right\Vert} \; \;
\left(\in L_{\sigma(j)}\right),
\endaligned
\tag 3.2
$$
$j=1,\dots,m$. Here $P_i$ denotes ($i=\sigma(j)$) the orthogonal projection
of $\Bbb R^d$ onto $L_i$. In this representation of the Euclidean path
$\gamma$ we are not at all bothered by the facts that

(a) the cylinders $a_j+C_{\sigma(j)}$ ($j=1,\dots,m$) may intersect each
other, or

(b) the path $\gamma$ itself may pass through certain cylinders 
$a_j+C_{\sigma(j)}$ without collision,

\noindent
because our investigation of Euclidean paths $\gamma$ (with a fixed symbolic
collision sequence $\Sigma$) will be a local, geometric analysis.

It is clear from 3.1 that the whole Euclidean path $\gamma=\gamma(\Sigma)$
is fully determined by the following data:

(i) the symbolic collision sequence 
$\Sigma=\left(\sigma(1),\dots,\sigma(m)\right)\in\{1,\dots,k\}^m$;

(ii) the translation vectors $a_j\in L_{\sigma(j)}$, $j=1,\dots,m$;

(iii) and by the initial (unit) velocity $V_0=\dot\gamma(0)\in\Bbb S^{d-1}$.

Therefore, the set $\Gamma=\Gamma(\Sigma)$ of all Euclidean paths 
$\gamma=\gamma(\Sigma)$ is naturally embedded into the product manifold
$\Bbb S^{d-1}\times \prod_{j=1}^{m}L_{\sigma(j)}$ as an open submanifold by
the mapping 

$$
\Psi:\, \Gamma\to\Bbb S^{d-1}\times \prod_{j=1}^{m}L_{\sigma(j)},
$$

$$
\Psi(\gamma)=\left(\dot\gamma(0);a_1(\gamma),\dots,a_m(\gamma)\right).
$$
In this way $\Gamma$ inherits a real analytic manifold structure from the 
ambient space $\Bbb S^{d-1}\times \prod_{j=1}^{m}L_{\sigma(j)}$. Set

$$
\aligned
\Gamma\left(\Sigma,\vec a\right)=\Gamma\left(\Sigma,a_1,\dots,a_m 
\right)= \\
\left\{\gamma\in\Gamma:\, \exists\; a\in\Bbb R^d\text{ such that }
a_j(\gamma)-a_j=P_{\sigma(j)}(a)\text{ for } j=1,\dots,m \right\},
\endaligned
\tag 3.3
$$
($a_j\in L_{\sigma(j)}$ are given) as the closed submanifold of $\Gamma$
corresponding to the given relative positions of the cylinders
$a_j+C_{\sigma(j)}$. It is easy to see that 
(if $\Gamma\left(\Sigma,\vec a\right)\ne\emptyset$) 
$\Gamma\left(\Sigma,\vec a\right)$ is a closed submanifold of 
$\Gamma(\Sigma)$ whose dimension is 
$2d-1-\text{dim}\left(\cap_{j=1}^m A_{\sigma(j)}\right)$.
{\it Throughout the paper we will only consider non-empty submanifolds
$\Gamma\left(\Sigma,\vec a\right)$.}

We need to introduce a special family of small perturbations of Euclidean
paths $\gamma\in\Gamma(\Sigma)$ corresponding to the pure spatial translations
of the initial phase point, see also Section 3 of [S-Sz(1998)], especially
formula (3.16) and its vicinity. Since we have now the convention 
$\gamma(0)=0$, instead of translating the initial position $\gamma(0)$,
we translate the cylinders $a_j+C_{\sigma(j)}$ by the same vector
$a\in\Bbb R^d$.

\subheading{\bf Definition 3.3-a} For $a\in\Bbb R^d$ ($||a||$ is small)
and $\gamma\in\Gamma=\Gamma(\Sigma)$ denote by $T_a(\gamma)=\delta$ the
uniquely defined element $\delta$ of $\Gamma$ for which
$\dot\delta(0)=\dot\gamma(0)$ and 
$a_j(\delta)-a_j(\gamma)=P_{\sigma(j)}(a)$, $j=1,\dots,m$.

In other terms, this means that we uniformly translate every scattering
cylinder $a_j+C_{\sigma(j)}$ of $\gamma$ by the same vector $a\in\Bbb R^d$,
which essentially amounts to the same thing as if we translated the initial
position by the vector $-a$.

In accordance with the part "Characterization of the Positive Subspace
of the Second Fundamental Form" in Section 3 of [S-Sz(1998)], we introduce
the following notions:

(a) the velocity process (history) $(V_0,V_1,\dots,V_m)$ of 
$\gamma\in\Gamma(\Sigma)$, where $V_0=\dot\gamma(0)$ and 
$V_j=\dot\gamma(t_j+0)$, $j=1,\dots,m$;

(b) the orthogonal reflection $h_j$ of $\Bbb R^d$ across the hyperplane

$$
H_j=\left(\dot\gamma(t_j+0)-\dot\gamma(t_j-0)\right)^\perp,
$$
$j=1,\dots,m$.

Note that the translated hyperplane $\gamma(t_j)+H_j$ is just the tangent
hyperplane of the boundary of the cylinder $a_j(\gamma)+C_{\sigma(j)}$ at
the point of reflection $\gamma(t_j)$. The collection of all possible
hyperplanes $H_j=H_j(\gamma)$ ($\supset A_{\sigma(j)}$) arising this way 
makes up the space $\Cal P_j$, being naturally diffeomorphic to the
($\nu_j-1$)-dimensional real projective space 
$\Bbb P^{\nu_j-1}(\Bbb R)$, where $\nu_j=\text{dim}L_{\sigma(j)}$, see
also Section 3 of [S-Sz(1998)]. Given an arbitrary sequence 
$(V_0;h_1,\dots,h_m)\in\Bbb S^{d-1}\times\prod_{j=1}^m \Cal P_j$,
one naturally defines the velocities 
$V_j=V_0\cdot h_1\cdot\dots\cdot h_j$ ($j=0,\dots,m$),
i. e. the image of $V_0$ under the composite
action $h_1\cdot\dots\cdot h_j$ of the reflections $h_1,\dots,h_j$.
(Here, by convention, the reflection $h_1$ is to be applied first.) Set

$$
\Phi\left(V_0;h_1,\dots,h_m\right)=V_m=
V_0\cdot h_1\cdot\dots\cdot h_m,
\tag 3.3-b
$$
cf. (3.17) of [S-Sz(1998)]. 

Let us observe that in the current representation of the Euclidean path
$\gamma\in\Gamma(\Sigma)$ with $\gamma(0)=0$, the notion of the neutral
space $\Cal N_0(\gamma)=\Cal N(\gamma)$ (cf. definition 2.1 in [S-Sz(1998)])
is redefined as follows:

$$
\Cal N(\gamma)=\left\{a\in\Bbb R^d:\, \exists\, \delta>0\text{ such that }
\forall\, \epsilon\in (-\delta,\delta) \;\; V_m\left(T_{\epsilon a}(\gamma)
\right)=V_m(\gamma)\right\}.
\tag 3.4
$$
For any vector $a\in\Bbb R^d$ and any Euclidean path $\gamma\in\Gamma(\Sigma)$
we introduce the following derivative:

$$
\aligned
\partial_a V_m=(\partial_a V_m)(\gamma)= \\
\lim_{\epsilon\to 0}\epsilon^{-1}\cdot\left(V_m\left(T_{\epsilon a}(\gamma)
\right)-V_m\left(\gamma\right)\right).
\endaligned
\tag 3.5
$$
In accordance with the notations of Proposition 3.18 of [S-Sz(1998)],
the subspace

$$
\Cal W_+=\Cal W_+(\gamma)=\left\{(\partial_aV_m)(\gamma):\, a\in\Bbb R^d
\right\}
\tag 3.6
$$
is precisely the positive subspace of the second fundamental form $W$ of the
image $S^t(B)$ ($t>t_m$) of the parallel "beam of light" 

$$
B=\left\{x=(q,v_0)\in\Bbb R^d\times\Bbb R^d:\, v_0=v_0(\gamma),\;
q\perp v_0,\; ||q||<\epsilon\right\}
$$
under the action $S^t(\, .\, )$ of the Euclidean cylindric billiard flow
determined by the cylinders $a_j(\gamma)+C_{\sigma(j)}$ generating the
collisions near time moments $t_j=t_j(\gamma)$, see also formula (3.16)
and the accompanying text in [S-Sz(1998)]. It is well known that the second
fundamental form $W$ is symmetric and positive semi-definite, thus we get

\subheading{\bf Proposition 3.7} The orthogonal complement
$\left(\Cal W_+(\gamma)\right)^\perp$ of $\Cal W_+(\gamma)$ is equal to the
image $\Cal N_0(\gamma)\cdot h_1\cdot\dots\cdot h_m$ of the neutral space
under the composite action $h_1\cdot\dots\cdot h_m$ of the reflections
$h_j=h_j(\gamma)$.

Besides the positive subspace $\Cal W_+(\gamma)$ we will need to use another
subspace of $\Bbb R^d$ associated with $\gamma$. Let us consider an arbitrary
vector $\vec b=(b_1,\dots,b_m)\in\prod_{j=1}^m L_{\sigma(j)}$. For any
Euclidean path $\gamma\in\Gamma(\Sigma)$ we denote by
$\delta=\Cal T_{\vec b}(\gamma)\in\Gamma(\Sigma)$ the uniquely determined
Euclidean path for which $V_0(\delta)=V_0(\gamma)$ and
$a_j(\delta)-a_j(\gamma)=b_j$, $j=1,\dots,m$. In other words, the perturbed
path $\delta$ corresponds to the translations of the cylinders
$a_j(\gamma)+C_{\sigma(j)}$ by the vectors $b_j\in L_{\sigma(j)}$. We note
here that --- since our analysis of Euclidean paths is local --- we are
only interested in {\it small} perturbations $\Cal T_{\vec b}$, so that
no problem arises concerning of the smoothness of the Euclidean cylindric
billiard flow. 

Set

$$
\aligned
\partial_{\vec b}V_m(\gamma)=\lim_{\epsilon\to 0}\epsilon^{-1}\cdot
\left[V_m\left(\Cal T_{\epsilon\vec b}(\gamma)\right)-V_m(\gamma)
\right], \\
\tilde{\Cal W}_+(\gamma)=\left\{\partial_{\vec b}V_m(\gamma):\,
\vec b\in\prod_{j=1}^m L_{\sigma(j)}\right\}.
\endaligned
\tag 3.8
$$
It is clear that $\Cal W_+(\gamma)\subset\tilde{\Cal W}_+(\gamma)$ and

$$
\tilde{\Cal W}_+(\gamma)=\text{Im}\left[\frac{\partial\Phi}
{\partial\tilde{\Cal P}}\left(V_0(\gamma);h_1(\gamma),\dots,h_m(\gamma)
\right)\right],
\tag 3.9
$$
where the right-hand-side of (3.9) denotes the image space of the partial
derivative of $\Phi:\, \Bbb S^{d-1}\times\tilde{\Cal P}\to\Bbb S^{d-1}$
with respect to the second factor 
$\tilde{\Cal P}=\prod_{j=1}^m \Cal P_j$, where the mapping
$\Phi(V_0;h_1,\dots,h_m)=V_0\cdot h_1\cdot\dots\cdot h_m$ is defined above,
see also (3.17) and the paragraph preceding Proposition 3.18 in [S-Sz(1998)].
The reason why the two sides of (3.9) coincide is that, by independently
translating the cylinders $a_j(\gamma)+C_{\sigma(j)}$ ($j=1,\dots,m$) 
one-by-one by the vectors $\epsilon\cdot b_j$, we can 
independently and arbitrarily perturb the reflections $h_j=h_j(\gamma)$,
as well. This argument immediately proves

\subheading{\bf Proposition 3.10} The mapping

$$
\Theta:\, \Gamma(\Sigma)\to \Bbb S^{d-1}\times\prod_{j=1}^m \Cal P_j,
$$
defined by 
$\Theta(\gamma)=\left(V_0(\gamma);h_1(\gamma),\dots,h_m(\gamma)\right)$
is a submersion (i. e. its derivative is surjective at every point) and,
hence, it is an open mapping. \qed

We cite here the fundamental assertion of Proposition 3.18 from
[S-Sz(1998)]:

\subheading{\bf Proposition 3.11} For every Euclidean path 
$\gamma\in\Gamma(\Sigma)$ the subspaces $\Cal W_+(\gamma)$ and
$\tilde{\Cal W}_+(\gamma)$ are equal. \qed

\subheading{\bf Remark} Observe that --- although Proposition 3.18 of
[S-Sz(1998)] was originally formulated and proven for paths of cylindric 
billiards in a torus, the entire proof obviously carries over to the
Euclidean case without any significant change.

Let us introduce now the following, useful notions of typical dimensions:

$$
\aligned
\Delta(\Sigma)=\max\Sb \gamma\in\Gamma(\Sigma)\endSb \dim\Cal W_+(\gamma)=
\max\Sb \gamma\in\Gamma(\Sigma)\endSb \dim\tilde{\Cal W}_+(\gamma)= \\
\max\left\{\text{dim}\text{Im}\left[\frac{\partial\Phi}
{\partial\tilde{\Cal P}}\left(V_0;h_1,\dots,h_m\right)\right]:\,
\left(V_0;h_1,\dots,h_m\right)\in\Bbb S^{d-1}\times\prod_{j=1}^m \Cal P_j
\right\},
\endaligned
\tag 3.12
$$

$$
\aligned
\Delta(\Sigma,\vec a)=\Delta(\Sigma;a_1,\dots,a_m)= \\
\max\left\{\text{dim}\Cal W_+(\gamma):\, \gamma\in\Gamma
(\Sigma;a_1,\dots,a_m)\right\}.
\endaligned
\tag 3.13
$$
For the definition of the non-empty, closed submanifold
$\Gamma(\Sigma;a_1,\dots,a_m)$, see also (3.3) above. We note that in 
the first equation of (3.12) we used Proposition 3.11, while in the second
equation of (3.12) we took advantage of (3.9) and Proposition 3.10.

The simple proof of the next proposition uses a quite common algebraic
argument.

\subheading{\bf Proposition 3.14} There exist three open sets with full
measure $\Cal O_1\subset\Gamma(\Sigma)$, 
$\Cal O_2\subset\Bbb S^{d-1}\times\tilde{\Cal P}$, and
$\Cal O_3\subset\Gamma(\Sigma;a_1,\dots,a_m)$ such that 

(i) $\text{dim}\Cal W_+(\gamma)=\Delta(\Sigma)$ for every 
$\gamma\in\Cal O_1$,

(ii)

$$
\text{dim}\text{Im}\left[\frac{\partial\Phi}{\partial\tilde{\Cal P}}
\left(V_0;h_1,\dots,h_m\right)\right]=\Delta(\Sigma)
$$
for every $\left(V_0;h_1,\dots,h_m\right)\in\Cal O_2$, and

(iii) $\text{dim}\Cal W_+(\gamma)=\Delta(\Sigma;a_1,\dots,a_m)$ for every
$\gamma\in\Cal O_3$.

\subheading{\bf Proof} We will only present here a brief sketch of the proof
for the first statement, for the arguments proving the other two are
analoguous.

The openness of $\Cal O_1\subset\Gamma(\Sigma)$ follows from the continuous
dependence of the linear generators $\partial_{e_i}V_m(\gamma)$ 
($i=1,\dots,d$; $e_i$ is the $i$-th standard unit vector in $\Bbb R^d$)
on $\gamma$, in other words, it follows from the lower semi-continuity
of the dimension function $\text{dim}\Cal W_+(\gamma)$.

The fact that the open set $\Cal O_1\subset\Gamma(\Sigma)$ has full measure
in $\Gamma(\Sigma)$ (more precisely: its complement is a countable union
of smooth, proper submanifolds of $\Gamma(\Sigma)$) follows from the following
observations: The coordinates of the linear generators
$\partial_{e_i}V_m(\gamma)$ ($i=1,\dots,d$) of the space $\Cal W_+(\gamma)$
are algebraic functions of the coordinates of

$$
\gamma=\left(V_0(\gamma);a_1(\gamma),\dots,a_m(\gamma)\right)\in
\Bbb S^{d-1}\times\prod_{j=1}^m L_{\sigma(j)}.
$$
These algebraic functions only contain constants, rational operations
(field operations), and square roots. The reason why this is indeed so comes
from the similar algebraic nature of the cylindric billiard dynamics: We are
dealing with circular cylinders as scatterers. Therefore, the kinetic data of
the process $\gamma$ itself (i. e. the time moments $t_j=t_j(\gamma)$,
the positions $\gamma(t_j)$, and the velocities $V_j=\dot\gamma(t_j+0)$) are
also algebraic functions of the above type of initial variables
$V_0(\gamma)$ and $a_j(\gamma)$. Recall that the time moment $t_j$ is
iteratively determined by the earlier kinetic variables as the smaller root
$\tau$ of the quadratic equation 

$$
\left\Vert P_{\sigma(j)}\left[\gamma(t_{j-1})+(\tau-t_{j-1})\dot
\gamma(t_{j-1}+0)-a_j\right]\right\Vert^2=r^2_{\sigma(j)},
\tag 3.15
$$
$j=1,\dots,m$. (Here we use the natural convention $t_0=0$.)
Note that the solutions of the equation (3.15) in the iterative process of
computing the variables $t_j$, $\gamma(t_j)$, and $\dot\gamma(t_{j}+0)$
($j=1,\dots,m$) is the only point where the square root enters the whole
process: all the other variables can be then expressed by rational operations.
For more details, see Section 3 of [S-Sz(1999)].

Consider now the $d\times d$ matrix

$$
M(\gamma)=\left(\partial_{e_1}V_m(\gamma),\dots,\partial_{e_d}V_m(\gamma)
\right)
$$
the entries of which are algebraic functions of the coordinates of the
variable

$$
\gamma=\left(V_0(\gamma);a_1(\gamma),\dots,a_m(\gamma)\right)\in\Bbb S^{d-1}
\times\prod_{j=1}^m L_{\sigma(j)}.
$$
The relation $\gamma\not\in\Cal O_1$ precisely means that 
$\text{rank}\left(M(\gamma)\right)<\Delta(\Sigma)$, i. e. every
$\Delta(\Sigma)\times\Delta(\Sigma)$ sized minor of $M(\gamma)$ is zero.
Since these minors are also algebraic functions of $\gamma$, and at least one
of them is not identically zero (because the value $\Delta(\Sigma)$ is
attained as $\text{rank}\left(M(\gamma)\right)$ for some 
$\gamma\in\Gamma(\Sigma)$), we get that the complement of $\Cal O_1$ in
$\Gamma(\Sigma)$ is indeed a countable union of proper, smooth submanifolds
of $\Gamma(\Sigma)$. (It as an algebraic set.) \qed

The next lemma effectively utilizes Proposition 3.14 and the theorem on
mappings with constant rank from the calculus of several variables.

\subheading{\bf Lemma 3.16} Let $\gamma\in\Cal O_1$
($\subset\Gamma(\Sigma)$), and a small number $\epsilon_0>0$ be given.
(We only study small perturbations.) Consider the following two sets of
final velocities $V_m$:

$$
\Cal V_1=\Cal V_1(\gamma,\Sigma,\epsilon_0)=\left\{V_m\left(T_a(\gamma)
\right):\, a\in\Bbb R^d,\; ||a||<\epsilon_0\right\},
$$

$$
\Cal V_2=\Cal V_2(\gamma,\Sigma,\epsilon_0)=\left\{V_m\left(T_{\vec b}
(\gamma)\right):\, \vec b=(b_1,\dots,b_m)\in\prod_{j=1}^m L_{\sigma(j)},
\; \max_{j}||b_j||<\epsilon_0\right\}.
$$
We claim that both $\Cal V_1$ and $\Cal V_2$ are 
$\Delta(\Sigma)$-dimensional, smooth manifolds containing $V_m(\gamma)$
(as an interior point),
and these manifolds coincide in a neighbourhood of the point $V_m(\gamma)$.

\subheading{\bf Proof} Both mappings

$$
a\longmapsto V_m\left(T_a(\gamma)\right)\quad (||a||<\epsilon_0)
$$
and

$$
\vec b\longmapsto V_m\left(T_{\vec b}(\gamma)\right)\quad
(\max_{j}||b_j||<\epsilon_0)
$$
have derivatives with constant rank $\Delta(\Sigma)$. Therefore, by the
mentioned theorem on mappings with constant rank (see, for instance, 
Theorem 15.5, Chapter I of [H(1978)]) the sets $\Cal V_1$ and $\Cal V_2$
are indeed $\Delta(\Sigma)$-dimensional, smooth, embedded submanifolds of 
$\Bbb R^d$ for small enough $\epsilon_0>0$. Since $\Cal V_1$ is obviously a
subset of $\Cal V_2$ in a neighbourhood of $V_m(\gamma)$ and these two
smooth manifolds have the same dimension, they must coincide in a
neighbourhood of the point $V_m(\gamma)$. \qed

The main result of this section is

\subheading{\bf Key Lemma 3.17} Assume that

$$
\vec a=(a_1,\dots,a_m)\in\prod_{j=1}^m L_{\sigma(j)}
$$
is such a multi-vector that $\Gamma(\Sigma,\vec a)\not=\emptyset$. Then the
typical dimensions of $\Cal W_+$ in $\Gamma(\Sigma)$ and 
$\Gamma(\Sigma,\vec a)$ are equal, i. e. 
$\Delta(\Sigma)=\Delta(\Sigma,\vec a)$.

\subheading{\bf Proof} Induction on the length $m$ of 
$\Sigma=\left(\sigma(1),\dots,\sigma(m)\right)$. For $m=1$ the assertion is
obviously true, for $\Gamma(\Sigma)=\Gamma(\Sigma,a_1)$. 

Assume now that $m>1$ and the key lemma has been proven for
$m'=1,\dots,m-1$. Consider and fix a symbolic sequence
$\Sigma=\left(\sigma(1),\dots,\sigma(m)\right)$ of length $m$ and a 
multi-vector $\vec a=(a_1,\dots,a_m)$ for which
$\Gamma(\Sigma,\vec a)\not=\emptyset$.

Denote by $\Sigma'$ the truncated sequence 
$\left(\sigma(1),\dots,\sigma(m-1)\right)$. Throughout the proof of the key
lemma, for $\gamma\in\Gamma(\Sigma)$ we denote by $\gamma'$ the following,
truncated Euclidean path: $\gamma'(t)=\gamma(t)$ for 
$0\le t\le t_{m-1}(\gamma)$, and
$\gamma'(t)=\gamma(t_{m-1})+(t-t_{m-1})\dot\gamma(t_{m-1}+0)$ for
$t\ge t_{m-1}(\gamma)$.

Select and fix an element $\gamma_0\in\Gamma(\Sigma,\vec a)$. By using the
induction hypothesis and the Fubini theorem, we can assume that the
truncated Euclidean path $\gamma'_0\in\Gamma(\Sigma';a_1,\dots,a_{m-1})$
belongs to the typical set $\Cal O_1(\Sigma')$ of $\Gamma(\Sigma')$
and, moreover, the following additional property also holds true:

$$
\cases
\text{For almost every selection of vectors }c_j\in L_{\sigma(j)}\quad
(j=1,\dots,m) \\
\text{the Euclidean path }
\delta=\left(V_0(\gamma);c_1,\dots,c_m\right)\in\Gamma(\Sigma)
\text{ (if exists!) belongs to the} \\
\text{typical set }\Cal O_1(\Sigma), \text{ and the truncated path }
\delta' \text{ is an element of } \Cal O_1(\Sigma'),
\endcases
\tag 3.18
$$
see Proposition 3.14 for the notion of the typical set $\Cal O_1$.

Select and fix a small number $\epsilon_1>0$. Its sufficient smallness will
be clarified later in the proof. There is now a perturbation
$\gamma_1\in\Gamma(\Sigma)$ of $\gamma_0$ with $V_0(\gamma_1)=V_0(\gamma_0)$
and $\Vert a_j(\gamma_1)-a_j(\gamma_0)\Vert<\epsilon_1$ ($j=1,\dots,m$) such
that $\gamma_1\in\Cal O_1(\Sigma)$ and $\gamma'_1\in\Cal O_1(\Sigma')$.
We note here that the relation $V_0(\gamma_1)=V_0(\gamma_0)$ can be achieved
just because of (3.18).

Consider and compare the two nearby Euclidean paths
$\gamma'_0,\, \gamma'_1\in\Cal O_1(\Sigma')$. Here
$\gamma'_0\in\Gamma(\Sigma';a_1,\dots,a_{m-1})$ also holds and, if the number
$\epsilon_1>0$ has been chosen small enough, the velocity
$V_{m-1}(\gamma'_1)=V_{m-1}(\gamma_1)$ belongs to the small open neighbourhood
$U_0\subset\Bbb R^d$ of the velocity $V_{m-1}(\gamma'_0)=V_{m-1}(\gamma_0)$
in which the sets $\Cal V_1=\Cal V_1\left(\gamma'_0,\Sigma',\epsilon_0\right)$
and $\Cal V_2=\Cal V_2\left(\gamma'_0,\Sigma',\epsilon_0\right)$ are 
$\Delta(\Sigma')$-dimensional, smooth manifolds and they coincide:
$\Cal V_1\cap U_0=\Cal V_2\cap U_0$, see Lemma 3.16. (Here we can see that
the number $\epsilon_0>0$ should be chosen first for $\gamma'_0$, according
to Lemma 3.16, and then $\epsilon_1>0$ must be selected small enough in order
to ensure the above properties.) Now we have

$$
V_{m-1}(\gamma_1)=V_{m-1}(\gamma'_1)\in\Cal V_1\cap U_0=\Cal V_2\cap U_0
\tag 3.19
$$
and, therefore, there exists a small perturbation

$$
\gamma_2=T_a(\gamma_0)\in\Gamma(\Sigma)\quad (||a||<\epsilon_0)
\tag 3.20
$$
for which $V_{m-1}(\gamma_2)=V_{m-1}(\gamma_1)$. If the first selected number
$\epsilon_0>0$ was chosen small enough then, necessarily, we have that
$\gamma'_2\in\Cal O_1(\Sigma')$, i. e. it is a typical Euclidean path for
$\Sigma'$. Consider now the three $\Sigma'$-typical Euclidean paths
$\gamma'_0,\, \gamma'_1,\, \gamma'_2\in\Cal O_1(\Sigma')$. Their neutral
linear spaces (measured now right after the collision $\sigma(m-1)$) are

$$
\Cal N(\gamma'_i)=\left(\Cal W_+(\gamma'_i)\right)^\perp,\quad (i=1,2,3),
$$
see Proposition 3.7. By the generic nature $\gamma'_i\in\Cal O_1(\Sigma')$
($i=1,2,3$) of $\gamma'_i$ we get

$$
\text{dim}\Cal W_+(\gamma'_i)=\Delta(\Sigma'), \quad i=1,2,3.
\tag 3.21
$$
On the other hand, the space $\Cal W_+(\gamma'_i)$ ($i=1,2,3$) is clearly
equal to the tangent space of the manifold $\Cal V_1\cap U_0=\Cal V_2\cap U_0$
at the point $V_{m-1}(\gamma'_i)$. Since 
$V_{m-1}(\gamma'_1)=V_{m-1}(\gamma'_2)$, we have that

$$
\cases
\Cal W_+(\gamma'_1)=\Cal W_+(\gamma'_2)\text{ and, therefore,} \\
\Cal N(\gamma'_1)=\Cal N(\gamma'_2)=\left(\Cal W_+(\gamma'_1)\right)^\perp.
\endcases
\tag 3.22
$$
The neutral space $\Cal N(\gamma_i)$ of the Euclidean path $\gamma_i$
($i=1,2$; the spaces $\Cal N(\gamma_i)$ are now measured between the
collisions $\sigma(m-1)$ and $\sigma(m)$) can be obtained obviously as the
intersection

$$
\Cal N(\gamma_i)=\Cal N(\gamma'_i)\cap\left(\Bbb R\cdot V_{m-1}(\gamma_i)+
A_{\sigma(m)}\right),
\tag 3.23
$$
$i=1,2$. Since the right-hand-sides of (3.23) are identical for $i=1$ and
$i=2$, we obtain that $\Cal N(\gamma_1)=\Cal N(\gamma_2)$ and, since
$\gamma_1\in\Cal O_1(\Sigma)$ is typical with respect to $\Sigma$, we have
that $\text{dim}\Cal N(\gamma_2)=\Delta(\Sigma)$, i. e.
$\gamma_2\in\Cal O_1(\Sigma)$. Taking into account (3.20), we see that
$\gamma_2\in\Gamma(\Sigma;\vec a)$, thus
$\Delta(\Sigma;\vec a)=\Delta(\Sigma)$, as claimed. The proof of Lemma 3.17
is now complete. \qed

\subheading{\bf Corollary 3.24} Suppose that $S^{[a,b]}x_0$ is a non-singular,
finite trajectory segment of the genuine, toroidal, cylindric billiard flow
$\flow$ with the collision sequence 
$\Sigma=\left(\sigma(1),\dots,\sigma(m)\right)$ for which

$$
\max\left\{\text{dim}\text{Im}\left[\frac{\partial\Phi}
{\partial\tilde{\Cal P}}\left(V_0;h_1,\dots,h_m\right)\right]:\,
\left(V_0;h_1,\dots,h_m\right)\in\Bbb S^{d-1}\times\prod_{j=1}^m \Cal P_j
\right\}=d-1,
\tag 3.25
$$
see also (3.12). (The numbers $a$ and $b$ are supposed to be non-collision
moments of time.) Then there is an open neighbourhood $U$ of $x_0$ in
$\bold M$ and there is a closed, proper (i. e. of codimension at least one)
algebraic set $F\subset U$ such that $S^{[a,b]}y$ is hyperbolic (sufficient)
for every $y\in U\setminus F$. \qed

\subheading{\bf Definition 3.26} Based upon the above corollary, we will say
that the symbolic sequence $\Sigma=\left(\sigma(1),\dots,\sigma(m)\right)$
is {\it combinatorially rich for one codimension} if (3.25) holds true.

\subheading{\bf Corollary 3.27} Theorem 5.1 of [S-Sz(1999)] along with
Corollary 3.24 imply that {\it every} hard ball system $\flow$ is
completely hyperbolic! Therefore, by Pesin's theory generalized to completely
hyperbolic dynamical systems with singularities [K-S(1986)], all ergodic
components of a hard ball system have positive measure, and the restriction
of the billiard flow to any ergodic component has the Bernoulli
property, see [C-H(1996)] and [O-W(1998)].
Thus, we see that the results of the present article are stronger than the
main theorem of [S-Sz(1999)] (where the complete hyperbolicity of
{\it almost every} hard ball system was proven), despite the fact that the
present approach does not use the rather involved algebraic machinery
of [S-Sz(1999)].

\bigskip \bigskip

\heading
4. Hyperbolicity Is Generic \\
(Proof of the theorem)
\endheading

\bigskip \bigskip

The goal of this section is to prove that in every {\it transverse}
cylindric billiard flow $\flow$ $\mu$-almost every phase point is hyperbolic
(in other words sufficient, see Section 2). This goal will be achieved through
the use of Corollary 3.24 by showing that the trajectory $S^{\Bbb R}x$ of
almost every phase point $x\in\bold M^0\cap\bold M^{\#}$ (for the definition
of the sets $\bold M^0$ and $\bold M^{\#}$ see Section 2 of this article)
contains infinitely many consecutive segments that are combinatorially rich
in the sense of Definition 3.26. It turns out, however, that in proving
this result the combinatorial richness described in 3.26 is not very
convenient for us, so we introduce the concept of a {\it transitive} 
(or, non-splitting) symbolic sequence $\Sigma$:

\subheading{\bf Definition 4.1} We say that the symbolic collision sequence
$\symb$ is {\it transitive} if the set of cylinders
$\left\{C_{\sigma(j)}:\, 1\le j\le m\right\}$ defines a transitive cylindric
billiard in the torus $\Bbb T^d=\Bbb R^d/\Cal L$ or, in other words, if the
system of base spaces $\left\{L_{\sigma(j)}:\, 1\le j\le m\right\}$ has the
Orthogonal Non-splitting Property, see 3.1--3.6 of [S-Sz(1998)], especially
3.3--3.4 and Theorem 3.6.

The next (elementary) lemma clarifies the relationship between the
transitivity of $\Sigma$ and its richness defined in 3.26.

\subheading{\bf Lemma 4.2} There exists an integer $C\in\Bbb N$ (depending
merely on $d$ and the transitive collection of base subspaces
$L_1,\dots,L_k\subset\Bbb R^d$) with the following properties: If a
symbolic sequence $\symb\in\{1,\dots,k\}^m$ contains at least $C$ 
consecutive, transitive subsequences, then the sequence $\Sigma$
is combinatorially rich as required by 3.26, i. e. formula (3.25)
holds true.

\subheading{\bf Proof} Let $\Cal T$ denote the set of all subsets
$T\subset\{1,2,\dots,k\}$ for which the collection of subspaces
$\{L_i:\, i\in T\}$ is transitive in $\Bbb R^d$. Let $|\Cal T|=n$
and $\Cal T=\{T_1,\dots,T_n\}$. For every $T_j$ ($1\le j\le n$) select
and fix a symbolic sequence 
$\Sigma^{(j)}=\left(\sigma^{(j)}(1),\dots,\sigma^{(j)}(m_j)\right)$
such that $\sigma^{(j)}(i)\in T_j$ ($i=1,\dots,m_j$), and $\Sigma^{(j)}$
is combinatorially rich in the sense of 3.26. Set
$C=n\cdot\max_{1\le j\le n}\{m_j\}$. If a symbolic sequence $\symb$
fulfills the condition of the lemma with the above constant $C$, then
there exists an index $j_0$ ($1\le j_0\le n$), and
$M=\max_{1\le j\le n}\{m_j\}$ consecutive subsegments
$\Sigma_1,\dots,\Sigma_M$ of $\Sigma$ with the property that the set of
labels in every $\Sigma_l$ ($1\le l\le M$) contains $T_{j_0}$ as a subset.
It is then clear that the rich sequence $\Sigma^{(j_0)}$ is a (rather
lacunary) subsequence of $\Sigma$, and being so, the considered symbolic
sequence $\symb$ is also combinarorially rich in the sense of 3.26. \qed

\bigskip

\heading
Eventually Splitting Trajectories
\endheading

\bigskip

\subheading{\bf Definition 4.3} We say that the positive semi-trajectory
$S^{(0,\infty)}x$ ($x\in\bold M^0$) splits according to the non-trivial
orthogonal splitting $\Bbb R^d=B_1\oplus B_2$ of $\Bbb R^d$ if for every
$t>0$ with $S^tx\in\partial C_i$ we have $L_i\subset B_1$ or $L_i\subset B_2$.

By keeping in mind Corollary 3.24 and Lemma 4.2, in order to prove our theorem
it is enough to obtain the following result, which is the analogue of 
Theorem 5.1 of [S-Sz(1999)].

\subheading{\bf Main Lemma 4.4} Assume that the cylindric billiard flow
$\flow$ has a transverse system $\{L_1,\dots,L_k\}$ of base spaces. Let
$\Bbb R^d=B_1\oplus B_2$ be a given non-trivial orthogonal splitting
of $\Bbb R^d$. We claim that the set

$$
S_{B_1,B_2}=\left\{x\in\bold M^0\cap\bold M^{\#}:\, S^{(0,\infty)}x
\text{ splits according to } B_1\oplus B_2 \right\}
$$
of phase points with $(B_1,B_2)$-splitting positive orbits has Liouville
measure zero, i. e. $\mu\left(S_{B_1,B_2}\right)=0$.
(Note that $\bold M^0$ denotes the set of all non-singular phase points,
while $\bold M^\#$ contains all phase points with infinitely many
non-tangential collisions in both time directions, see also Section 2.)

\subheading{\bf Proof} The rest of this section will be devoted to the proof
of the main lemma. The proof will be subdivided into a few lemmas.

Consider and fix an arbitrary phase point
$x_0\in S_{B_1,B_2}\setminus\partial\bold M$ 
($\subset\bold M^0\cap\bold M^\#$). We want to show that $x_0$ has an open
neighbourhood $U\subset\bold M\setminus\partial\bold M$ for which
$\mu\left(S_{B_1,B_2}\cap U\right)=0$.

We set

$$
I=\left\{i:\; 1\le i\le k,\; \exists\, t>0 \text{ such that }
S^tx_0\in\partial C_i \right\}.
$$
Plainly, $I\not=\emptyset$. By switching from $x_0$ to an image $S^tx_0$
($t>0$) if necessary, we can assume that for everi $i\in I$ there is an 
infinite sequence $t_n\nearrow\infty$ such that 
$S^{t_n}x_0\in\partial C_i$ ($n\in\Bbb N$), i. e. the set $I$ is already
stable.

Clearly, the Euclidean space $\Bbb R^d$ uniquely splits into an orthogonal
direct sum

$$
\Bbb R^d=\bigoplus_{j=1}^p E_j\oplus E_0,
\tag 4.5
$$
where

(i) for $j=1,\dots,p$ $\text{dim}E_j\ge 2$, and the base spaces
$\left\{L_i:\, i\in I,\; L_i\subset E_j\right\}$ enjoy the transitivity
(or the Orthogonal Non-splitting Property, see the definition right before
Lemma 3.3 in [S-Sz(1998)]) in $E_j$;

(ii) $\forall i\in I$ $\exists j$ ($1\le j\le p$) such that 
$L_i\subset E_j$.

Since the system $\{L_i:\, i\in I\}$ splits, by the assumed transverseness
of the entire system $\{L_1,\dots,L_k\}$ we have that there exists an index
$j_0\in\{1,\dots,k\}$ with the following property:

$$
P_{E^+}(A_{j_0})=E^+,
\tag 4.6
$$
where

$$
E^+=\bigoplus_{j=1}^p E_j=E_0^\perp,
$$
and $P_{E^+}$ denotes the orthogonal projection of $\Bbb R^d$ onto $E^+$.
Since $\text{dim}A_{j_0}\le d-2$, as a consequence, we get that

$$
\text{dim}E_0\ge 2.
\tag 4.7
$$

\subheading{\bf Remark 4.8} It follows easily from (i)--(ii) above that 
$p\ge1$, and the linear span $\text{span}\{L_i:\, i\in I\}$ is equal to the
space $E^+$, see also Remark 3.5 in [S-Sz(1998)]. As far as the special
index $j_0$ (featuring (4.6)) is concerned, we certainly have that 
$j_0\not\in I$, otherwise the projection $P_{E^+}(A_{j_0})$ would be 
orthogonal to the space $L_{j_0}\subset E^+$.

\subheading{\bf Definition 4.9} The $I$-dynamics $S_I^ty$ ($y\in\bold M$,
$t>0$) is defined as follows: $S_I^ty$ evolves according to the sub-billiard
system $\{C_i:\, i\in I\}$ in $\Bbb T^d=\Bbb R^d/\Cal L$, i. e. for $t>0$
we no longer remove the cylinders $\{C_i:\, i\not\in I\}$ from the
configuration space (i. e. we no longer considering them as scatterers)
but, instead, we allow
for the moving point $q\left(S_I^ty\right)$ to freely pass through the
transparent cylinders $C_i$ with $i\not\in I$. As to the notion of
sub-billiards, see subsection 2.1.

Obviously, in order to prove Main Lemma 4.4 it is enough to prove

\subheading{\bf Proposition 4.10} There exists an open neighbourhood
$U\subset\bold M\setminus\partial\bold M$ of $x_0$ in $\bold M$ such that

$$
\mu\left(\left\{y\in U\cap\bold M^0\cap\bold M^\# :\, \forall\; t>0
\quad S^ty=S_I^ty\right\}\right)=0.
$$

In the sequel we will just prove Proposition 4.10.

Fix the values of the partial kinetic energies 
$\varepsilon_j=\frac{1}{2}\Vert P_{E_j}(v)\Vert^2$ ($j=1,\dots,p$) and the
velocity $P_{E_0}(v)$. Introduce the notation $\jflow$ for the sub-billiard
flow determined by the index set 
$I_j=\left\{i:\, i\in I,\; L_i\subset E_j\right\}$ and by the given
kinetic energy $\varepsilon_j$, $j=1,\dots,p$. The configuration space of
this sub-billiard flow is naturally the torus $E_j/P_{E_j}(\Cal L)$
minus the intersections of the cylinders $\{C_i:\, i\in I_j\}$
with this torus, see also subsection 2.1. We note that the space $E_j$
here corresponds to the notation $E^+$ of 2.1.

Introduce also the notations $\Cal A_j=E_j^\perp \subset\Bbb R^d$,
$\Cal T_j=\Cal A_j/(\Cal A_j\cap\Cal L)$, and
$\Cal T_0=\bigcap_{j=1}^p \Cal T_j$ for $j=1,\dots,p$. Note that ---
as it is easy to see --- the subspaces $\Cal A_j$ are lattice subspaces,
thus $\Cal T_j$ are subtori of $\Bbb T^d=\Bbb R^d/\Cal L$, and $\Cal T_0$
is a closed subgroup of $\Bbb R^d/\Cal L$ being a finite extension of the
subtorus $E_0/(E_0\cap\Cal L)$.

\subheading{\bf Lemma 4.11} After fixing the values of the partial kinetic
energies \newline
$\varepsilon_j=\frac{1}{2}\Vert P_{E_j}(v)\Vert^2$ ($j=1,\dots,p$)
and the velocity $P_{E_0}(v)$, there exists a natural homomorphism of
dynamical systems

$$
\Psi:\; \left(\bold{M}_I,\{S_I^t\},\mu_I\right)\longrightarrow
\prod_{j=1}^p\jflow
$$
for which

(i) $\Psi$ is surjective;

(ii) two phase points $(q_1,v_1),\, (q_2,v_2)\in\bold M_I$ are mapped to the
same element by $\Psi$ if and only if $v_1=v_2$ and 
$q_1-q_2\in\Cal T_0$.

Therefore, the dynamical system $\left(\bold{M}_I,\{S_I^t\},\mu_I\right)$
is locally isomorphic to the direct product

$$
\prod_{j=1}^p\jflow
$$
multiplied by the uniform (conditionally periodic) motion in the torus
$E_0/(E_0\cap\Cal L)$ with the given velocity $P_{E_0}(v)$.

\subheading{\bf Proof} According to subsection 2.1, the sub-billiard flow
$\jflow$ is naturally a factor of $\left(\bold{M}_I,\{S_I^t\},\mu_I\right)$.
Denote by $\Psi_j$ the natural projection of the latter dynamical system
onto the former one, $j=1,\dots,p$. Thanks to the orthogonality of the
bases of cylinders in $\bold M_{j_1}$ and $\bold M_{j_2}$ ($j_1\ne j_2$),
we see that the $j_1$-part and $j_2$-part of the $S_I$-evolving phase
point $S_I^t y_0=y_t=(q_t,v_t)$ evolve independently. This shows that
the mapping $\Psi=(\Psi_1,\dots,\Psi_p)$ with the components $\Psi_j$
is a homomorphism between the dynamical systems 
$\left(\bold{M}_I,\{S_I^t\},\mu_I\right)$ and $\prod_{j=1}^p\jflow$.
It is obvious that the mapping $\Psi$ is surjective.

The only outstanding question is (ii) in the lemma. Assume, therefore, that
$\Psi(q_1,v_1)=\Psi(q_2,v_2)$. Since $P_{E_j}(v_1)=P_{E_j}(v_2)$ for
$j=0,1,\dots,p$, we immediately have that $v_1=v_2$. On the other hand,
the equation of the $q$-components of $\Psi_j(q_1,v_1)$ and 
$\Psi_j(q_2,v_2)$ precisely means that 
$q_1-q_2\in\Cal A_j/(\Cal A_j\cap\Cal L)=\Cal T_j$, $j=1,\dots,p$, i. e.
$q_1-q_2\in\Cal T_0=\bigcap_{j=1}^p \Cal T_j$. \qed

\subheading{Proof of Proposition 4.10} First of all, it is enough to prove
4.10 for fixed values of
$\varepsilon_j=\frac{1}{2}\Vert P_{E_j}(v)\Vert^2$ ($j=1,\dots,p$)
and the velocity $P_{E_0}(v)=v_0$. Thus, let us fix these values and prove
4.10 for the corresponding layer of the phase space.

Since for every $i\in I$ there is an infinite sequence 
$t_n\nearrow\infty$ such that $S^{t_n}x_0\in\partial C_i$, by applying
Lemma 4.2, property (i) after (4.5) and Corollary 3.24 for the sub-billiard
factor $\jflow$ ($j=1,\dots,p$), we obtain that there exists an open
neighbourhood $U\subset\bold M\setminus\partial\bold M$ of $x_0$ in
$\bold M$ and a proper, smooth submanifold $N\subset U$ such that 

$$
\cases
\text{the sub-billiard semi-orbit } \left\{S_j^ty:\; t>0\right\} \\
\text{is hyperbolic for every } y\in U\setminus N
\text{ and } j=1,\dots,p.
\endcases
\tag 4.12
$$
We note here that --- although the $I$-dynamics $S_I^ty$ ($y\in U$, $t>0$)
is not isomorphic to the direct product

$$
\prod_{j=0}^p \jflow
$$
(where $\left(\bold{M}_0,\{S_0^t\},\mu_0\right)$ is the conditionally
periodic motion in the torus $E_0/(E_0\cap\Cal L)$), but they are still
locally isomorphic according to Lemma 4.11. Therefore, in the small
neighbourhood $U$ the semi-orbit $S_I^ty$ can be written as

$$
S_I^ty=\left(S_0^ty_0,S_1^ty_1,\dots,S_p^ty_p\right)
\tag 4.13
$$
($y\in U$, $t>0$, $S_j^ty_j\in\bold M_j$) by using a local isomorphism
provided by Lemma 4.11. Thanks to (4.12) and the generalized Pesin theory
for hyperbolic dynamical systems with singularities [K-S(1986)], for
$\mu$-almost every phase point $y\in U\setminus N$ and $j=1,\dots,p$
the above component
$y_j$ of $y$ belongs to an ergodic component $C^{(j)}_{\alpha_j(y)}$
of the flow $\{S_j^t\}$ with the following properties:

$$
\mu_j\left(C^{(j)}_{\alpha_j(y)}\right)>0,
\tag 4.14
$$

$$
S_j^t|C^{(j)}_{\alpha_j(y)}\text{ is a mixing flow}.
\tag 4.15
$$
By considering generic phase points $y\in U\setminus N$, we can assume that
the fixed velocity $v_0=P_{E_0}(v)$ of the uniform motion $S_0^ty_0$ is
ergodic. Let us, therefore, denote by

$$
U(v_0,\varepsilon_1,\dots,\varepsilon_p,\alpha_1,\dots,\alpha_p)=
U(v_0,\vec\varepsilon,\vec\alpha)
$$
the set of all phase points $y=(q,v)\in(U\setminus N)\cap\bold M^0$ for
which $P_{E_0}(v)=v_0$, $\frac{1}{2}\Vert P_{E_j}(v)\Vert^2=\varepsilon_j$,
$\alpha_j(y)=\alpha_j$ for $j=1,\dots,p$, and the $I$-dynamics $S_I^ty$
is non-singular (just as $S^ty$) for $t>0$. We want to prove that

$$
\mu\left(\left\{y\in U(v_0,\vec\varepsilon,\vec\alpha):\, \forall\; t>0
\quad S^ty=S_I^ty \right\}\right)=0.
\tag 4.16
$$
The direct product flow

$$
\left(\bold{M}_0,\{S_0^t\},\mu_0\right)\times\prod_{j=1}^p
\left(C^{(j)}_{\alpha_j},\{S_j^t\},\mu_j|C^{(j)}_{\alpha_j}\right)
\tag 4.17
$$
(which governs the time evolution of $S_I^ty$, 
$y\in U(v_0,\vec\varepsilon,\vec\alpha)$, $S_0^t(y_0)=y_0+tv_0$) is ergodic
--- being the product of $p$ mixing flows and an ergodic one. The condition
$S^ty=S_I^ty$ ($\forall\; t>0$) specially means that the interior of the
cylinder $C_{j_0}$ (see (4.6)) is avoided. This is just the well studied
phenomenon of open set (ball) avoiding! The geometric condition (4.6) means
that for any given $p$-tuple of positions

$$
(q_1,\dots,q_p)\in\prod_{j=1}^p \left(E_j/P_{E_j}(\Cal L)\right)\cong
\prod_{j=1}^p \Bbb R^d/(\Cal L+\Cal A_j)
$$
one can find an element $\tilde q\in\Bbb R^d/\Cal L$ such that
$\Psi_j(\tilde q)=q_j$ ($j=1,\dots,p$) or, in other words, the natural
projection 

$$
\pi_j:\; \Bbb R^d/\Cal L\longrightarrow \Bbb R^d/(\Cal L+\Cal A_j)
\cong E_j/P_{E_j}(\Cal L)
$$
(see subsection 2.1) maps $\tilde q$ onto $q_j$, $\pi_j(\tilde q)=q_j$.
More precisely, the geometric condition (4.6) implies that for every element
$q\in\Bbb T^d=\Bbb R^d/\Cal L$ there exists another element 
$\tilde q\in\Bbb R^d/\Cal L$ for which 

$$
\tilde q\in(A_{j_0}/\Cal L)+t_{j_0}\subset\text{int}C_{j_0},
$$
and $\tilde q-q\in E_0/(E_0\cap\Cal L)$, i. e. even the actual connected
component of the inverse image $\Psi^{-1}\left((q_1,\dots,q_p)\right)$
(to contain $\tilde q$) can be specified arbitrarily.
(Recall that the translated subtorus $(A_{j_0}/\Cal L)+t_{j_0}$ is just
the axis of the cylinder $C_{j_0}$, see also the introduction.)
Especially, the phase space

$$
\left(E_0/(E_0\cap\Cal L)\right)\times\prod_{j=1}^p
C^{(j)}_{\alpha_j}
$$
of the flow in (4.17) has an intersection of positive measure with the 
interior of the "forbidden" cylinder $C_{j_0}$. Therefore, due to the
ergodicity of the product in (4.17), the event 
$\forall\, t>0 \quad S^ty=S_I^ty$
($y\in U(v_0,\vec\varepsilon,\vec\alpha)$) has indeed zero measure with
respect to the product measure in (4.17), consequently (4.16) is true.

This finishes the proof of Proposition 4.10 and Main Lemma 4.4. \qed

On the other hand, Main Lemma 4.4 together with Corollary 3.24 yield a
proof for the theorem of this article. \qed

\subheading{Corollary 4.18} It follows from the generalized Pesin theory
for hyperbolic dynamical systems with singularities [K-S(1986)] that
every transverse cylindric billiard system has at most countably many
ergodic components $C_{\alpha}$ (with positive measure), and the restrictions
$S^t|C_{\alpha}$ of the flow have the Bernoulli property, see [C-H(1996)]
and [O-W(1998)]

\bigskip

\subheading{\bf Concluding remark} The property of transverseness somehow
means that (in rough terms) the generator spaces $A_i$ of the cylinders
are big, as opposed to the condition ($A_i\cap A_j=\{0\}$ for $i\not=j$)
that was assumed by P. B\'alint in his Theorem 2.4 of [B(1999)].
Thus, we can say that --- in some sense --- the result of this article
is sort of complementary to B\'alint's Theorem 2.4. Out of these two result
it is the present one that applies to hard ball systems.

\enddocument